
\magnification=\magstep1
\baselineskip=15pt
\parskip=4pt

\def\C{{\bf C}}
\def\R{{\bf R}}
\def\l{\lambda}
\def\L{\Lambda}
\def\G{{\Gamma}}
\def\Re{{\rm Re}}
\def\Im{{\rm Im}}
\def\ind{{\rm ind}\,}

\font\auth=cmcsc10
\def\jour{\sl}
\def\vol{\bf}

\def\br{\hfill\break}

\topglue 1truein
\centerline{\bf
ANALYTIC CONTINUATION FROM A FAMILY OF LINES}
\bigskip
\centerline{{\auth Alexander Tumanov}}
\bigskip
\bigskip

{\narrower\narrower
{\bf Abstract.}
Given a function $f$ in the exterior of a convex curve
in the real plane, we prove that if the restrictions of $f$
to the tangent lines to the curve extend as entire functions,
then the function $f$ is an entire function of two
variables.

}
\bigskip\bigskip

We prove the following
\medskip
{\bf Theorem 1.
\sl
Let $\L\subset \R^2$ be a $C^2$-smooth convex curve with strictly
positive curvature.
Let $f$ be a complex function in the exterior of $\L$.
Denote by $l_\l$ the tangent line to $\L$ at $\l\in\L$.
Suppose for every $\l\in\L$, there is an entire function
$f_\l$ on $\C=\R^2$ such that the restrictions of $f$
and $f_\l$ to $l_\l$ coincide, that is,
$f|_{l_\l}=f_\l|_{l_\l}$.
Suppose that the map $(z,\l)\mapsto f_\l(z)$ is
continuous.
Then $f$ extends as an entire function on $\C^2$.
}
\medskip

The condition that $\L$ has positive curvature is
not essential. We add it for simplicity and convenience
of presentation. Instead of assuming that the map
$(z,\l)\mapsto f_\l(z)$ is continuous, we can
assume that it is locally bounded.
If $\L$ is real-analytic, then this condition can be
dropped (Theorem 6).
The author does not know whether it can be dropped
if $\L$ is merely smooth.

In the case $\L$ is the unit circle, the result
was obtained by Aguilar, Ehrenpreis and Kuchment [AEK]
as the characterization of the range of a version
of the Radon transform.
\"Oktem [\"O] gave another proof using a separate
analyticity result of Siciak [S2].
Ehrenpreis informed the author that he was able to
prove Theorem 1 for algebraic $\L$.
The proofs in [AEK, \"O]  don't go through for an arbitrary
convex curve $\L$. Our method is similar to that of [T].
It is based on Lewy's [L] proof of the classical
extension theorem of Kneser [K] and Lewy [L].

The author thanks Peter Kuchment for formulating
the problem and useful discussions.
The author thanks Guiseppe Zampieri for useful
discussions.
\bigskip
\centerline{*\quad*\quad*}
\medskip

Turning to the proof, we introduce
$$
\Sigma=\{ (z,w)\in\C^2: w=\bar z \}
$$
$$
h_\l=\{(z,\bar z): z\in l_\l \}\subset\Sigma
$$
Let $L_\l$ be the complex line in $\C^2$ containing
the real line $h_\l$.
Define a hypersurface
$$
M=\bigcup_{\l\in\L}L_\l.
$$
Define a function $F$ on $M$ that consists of all the
extensions $f_\l$ of the original function $f$,
that is,
$$
F(z,w)=f_\l(z) \quad{\rm if}\quad (z,w)\in L_\l.
$$
We will show that $F$ extends as an entire function
on $\C^2$. Then the original function $f$ extends as
an entire function of $z$ and $\bar z$, hence, as an entire
function of $\Re z$ and $\Im z$ as desired.

We first describe the geometry of $M$.
Let $\L_i$ and $\L_e$ denote the two connected components
of $\C \setminus \L$, the interior and exterior of $\L$.
Suppose $0\in\L_i$.
Likewise, we put
$$
\Pi=\{ (z,w)\in\Sigma: z\in\L \},\qquad
\Pi_e=\{ (z,w)\in\Sigma: z\in\L_e \}.
$$
Since $\L$ has positive curvature, the line $L_\l$
moves transversally to itself as $\l$ traces $\L$
except at $\Pi$, where the $L_\l$ is sliding
along itself. Then $M\setminus\Pi$ is a smooth
immersed hypersurface. The lines $L_\l$ intersect
only on $\Pi_e$, and every point of $\Pi_e$ is the point of
intersection of exactly two lines $L_\l$.
Hence, $M\setminus\Pi$ is a real Levi-flat hypersurface in
$\C^2$ with transverse self-intersection on $\Pi_e$.

We parametrize $\L$ by the arc length $s\mapsto \l\in\L$ in
counterclockwise direction, then $\l'=d\l/ds$ is the
unit tangent vector at $\l$, in particular $|\l'|=1$.
The tangent line $l_\l$ has the parametric equation
$$
z=\l+t\l',\quad
t\in\R.
$$
Then $h_\l\subset\Sigma$ has the parametric equations
$$
z=\l+t\l', \quad
w=\bar\l+t\bar\l',\quad
t\in\R.
$$
Eliminating $t$ we get
$$
{z-\l \over\l'}={w-\bar\l\over\bar\l'}\in\R.
$$
Then the complexification $L_\l$ of $h_\l$ has the equation
$$
{z-\l \over\l'}={w-\bar\l\over\bar\l'}.
$$
Let $l_\l^\pm$ be the two components of $\C\setminus l_\l$,
with $0\in l_\l^+$, that is,
$$
l_\l^\pm=\{ z\in\C: \pm\Im{z-\l\over \l'}>0 \}.
$$
We put
$$
L_\l^\pm=\{(z,w)\in L_\l: z\in l_\l^\pm  \},\qquad
M^\pm=\bigcup_{\l\in\L}L_\l^\pm.
$$
Then $M\setminus\Sigma=M^+\cup M^-$.
We study the geometry of $M$ by looking at the intersections
with coordinate lines
$$
\G_z=\{w\in\C: (z,w)\in M\},\qquad
\G_z^\pm=\{w\in\C: (z,w)\in M^\pm\},
$$
where $z\in\C$ is fixed.
The curve $\G_z$ is parametrized as
$$
\L\ni\l\mapsto w(\l)=\bar\l+(\bar\l')^2(z-\l).
$$
Note that for $\l_1\ne\l_2$,
$L_{\l_1}\cap L_{\l_2}=h_{\l_1}\cap h_{\l_2}
\subset\Sigma$. Then $\G_z$ can have
self-intersection only at the point $w=\bar z$.
The point $(z,\bar z)\in h_\l$ if and only if
${z-\l\over\l'}\in\R$, that is, $z\in l_\l$.

If $z\in\L_i$, then for every $\l\in\L$, we have
$\Im{z-\l\over \l'}>0$, hence $\G_z=\G_z^+$
is a simple closed curve and $\G_z^-=\emptyset$.

If $z\in\L_e$, then there are two tangent lines
$l_\l$ through $z$.
The points of tangency split $\L$ into two arcs
$$
\L^\pm(z)=\{\l\in\L: w(\l)\in \G_z^\pm \}
=\{\l\in\L:\pm\Im{z-\l\over \l'}>0\}.
$$
The arc $\L^-(z)$ is the one closer to $z$.
Hence, the curve $\G_z$ consists of two simple
disjoint loops $\G_z^\pm$ initiating and
terminating at $\bar z$.

We claim that $\G_z^-$ lies inside $\G_z^+$.
Indeed, we have
$$
|w(\l)-(\bar\l')^2 z|=|\bar\l-\l(\bar\l')^2|\le
2|\l|\le 2m,\quad
\hbox{where} \quad
m=\max_{\l\in\L}|\l|.
$$
Denote by $\ind \phi$ the index or winding number
of a map $\phi:\L\to\C$ around 0 as $\L$ is traced
counterclockwise.
Then $\ind \l'=1$, $\ind (\bar\l')^2 z=-2$.
In view of the above estimate of $w(\l)$, for big $|z|$,
we have $\ind \G_z=-2$. This implies that
$\G_z^\pm$ lie inside one another.
The curve $\G_z^-$ contracts into $\bar z_0$
as $\L_e\ni z\to z_0\in\L$.  Hence, $\G_z^-$ is the
inner one.

Let $D_z^+$, $D_z^-$, and $D_z^0$ be the three connected
components of $\C\setminus \G_z$:
the exterior of $\G_z^+$, the interior of $\G_z^-$,
and the domain between them. If $z\in\L_i$,
then $D_z^-=\emptyset$.
Then
$$
\Omega^\nu=\bigcup_{z\in\C}\{z\}\times D_z^\nu,\quad
\nu=+, -, 0,
$$
are the three connected components
of $\C^2\setminus M$.
The estimate of $w(\l)$ implies
$
D_z^\pm\supset \{w\in\C: \pm(|w|-|z|)>2m \}
$,
hence,
$$
\Omega^\pm\supset \{(z,w)\in\C^2: \pm(|w|-|z|)>2m \}.
$$

By the hypotheses of Theorem 1, the function $F$
defined above is a continuous CR function on $M$.
\medskip

{\bf Lemma 2.}
{\sl
$F$ holomorphically extends into $\Omega^-$.
}

{\it Proof.}
We follow the proof of the H. Lewy [L] extension theorem.
We define the extension $\tilde F$ by the Cauchy
type integral
$$
\tilde F(z,w)={1\over 2\pi i}
\int_{\G_z^-} {F(z,\zeta)\,d\zeta\over \zeta-w},\quad
z\in\L_e,\quad
w\in D_z^-.
$$
Obviously, $\tilde F$ is holomorphic in $w$.
To see that $\tilde F$ is an extension of $F$,
we prove the vanishing of the moments
$$
m_n(z)=\int_{\G_z^-} \zeta^n F(z,\zeta)\,d\zeta, \quad
z\in\L_e,\quad
n\ge0.
$$
We first show that both $\tilde F$ and $m_n$ are
holomorphic in $z$.
Put $\Phi(z)=\int_{\G_z^-}H(z,\zeta)\,d\zeta$,
where $H$ is the integrand in the formulas for
$\tilde F$ and $m_n$.
By the Morera theorem, it suffices to show
that $\int_\gamma\Phi(z)\,dz=0$
for every small loop $\gamma$ in $\L_e$.

The torus
$T=\bigcup_{z\in\gamma}\{z\}\times \G_z^-\subset M$
bounds a solid torus
$S\subset M$ obtained by filling the loop
$\gamma$.
Then by Stokes' formula
$$
\int_\gamma\Phi(z)\,dz=
\int_T H(z,\zeta)\,d\zeta\wedge dz=
\int_S dH(z,\zeta)\wedge d\zeta\wedge dz=0
$$
because $H$ is a CR function on $M$.

Recall that $\G_z^-$ contracts into $\bar z_0$
as $\L_e\ni z\to z_0\in\L$.
Then $m_n(z)\to 0$ as $z\to z_0$.
By the boundary uniqueness theorem,
$m_n(z)= 0$ identically.
By the moment condition, $\tilde F$ is the
holomorphic extension of $F$ as desired.
The lemma is proved.
\medskip

Note that the interchange $z\leftrightarrow w$
leads to the interchange
$\Omega^+\leftrightarrow\Omega^-$.
Then $F$ also extends into $\Omega^+$.
It remains to prove
\medskip

{\bf Lemma 3.
\sl
$\tilde F$ holomorphically extends into
$\Omega^0$.

\it Proof.}
Let
$$
X_c=\{(z,w): zw=c\},
$$
where $c\in\C$, and $|c|$ is large.
The intersection
$K_c=X_c\cap M$ is (the image of)
the map
$K_c:\L\ni\l\mapsto X_c\cap L_\l$.
The curve $K_c$ can have self-intersection
only if $w=\bar z$, hence $c=|z|^2>0$.

We use $z$ as a coordinate on $X_c$, which
identifies $X_c$ with the punctured
complex plane $\C\setminus\{0\}$.
Solving the equations for $X_c\cap L_\l$, we find
$$
z={-ib\pm\sqrt{c-b^2}\over\bar\l'}, \quad
b=\Im(\l'\bar\l).
$$
If $c\notin\R$, then $K_c$ consists of two
disjoint loops corresponding to the two values
of the square root. They are different because
$K_c$ can't have self-intersections.

Since $|b|\le m$, then
both components of $K_c$ are within fixed distance
from the circle $|z|=\sqrt{|c|}$. Hence
they lie inside one another.
Let $X_c^+$, $X_c^-$, and $X_c^0$ be the three
components of $X_c\setminus K_c$, with
$X_c^+$ being the inmost component,
$X_c^-$ being the outmost component,
and $X_c^0$ being the annulus between them.
Then $X_c^\pm\subset\Omega^\pm$ because
$\Omega^\pm\supset \{(z,w)\in\C: \pm(|w|-|z|)>2m \}.$
Hence $X_c^0\subset \Omega^0$.

If $c\in\R^+$, the positive reals, then the two
curves corresponding to the two values of the square
root coincide and have the surprisingly simple
form
$$
K_c=\{(z,w): |z|=\sqrt{c}, w=\bar z \}.
$$
We now use the same argument as in Lemma 2 in
the coordinate system $(z,c=zw)$.
We define the extension $\tilde F$ to $X_c^0$
by the Cauchy integral along $K_c$.
The vanishing of the moments will take place
because the integrals along the two components
of $K_c$ will cancel in the limit as
$\C\setminus\R^+\ni c\to c_0\in\R^+$.
This will yield the holomorphic extension of
$F$ to $X_c$ for large $c$.
Then by the Hartogs extension theorem,
$\tilde F$ extends to the whole space $\C^2$.

Lemma 3 and Theorem 1 are now proved.

\medskip
In conclusion, we show how to get rid of
the continuity hypothesis in Theorem 1 if
$\L$ is real-analytic.
We use a separate analyticity result by Siciak [S1].

\medskip
{\bf Lemma 4.} (Siciak [S1])
{\sl
Let $E_\nu\subset\C \;\; (\nu=1,2)$ be the domain
bounded by the ellipse with foci $\pm r_\nu$, and let
$I_\nu=[-r_\nu, r_\nu]\subset\R$.
Let $X=(E_1\times I_2)\cup (I_1\times E_2)\subset\C^2$.
Let $f:X\to\C$ be separately holomorphic, that is,
for every $x_1\in I_1$, $f(x_1,\cdot)$ is holomorphic in $E_2$, and
for every $x_2\in I_2$, $f(\cdot,x_2)$ is holomorphic in $E_1$.
Then $f$ holomorphically extends into a neighborhood of $X$.
}
\medskip
Bernstein (1912) proved Lemma 4 under an additional boundedness
assumption. He gave a precise description of the domain
of the extension, but we don't need it here.
Lemma 5 implies a propagation of analyticity
result for separately holomorphic functions.
\medskip
{\bf Proposition 5.}
{\sl
Let $M$ be a real-analytic Levi-flat hypersurface
in $\C^2$. Let $S\subset M$ be a connected complex
curve. Let $F:M\to\C$ be a separately holomorphic
function, that is, $F$ is holomorphic on complex
curves in $M$. Suppose $F$
holomorphically extends to a neighborhood
of $p\in S$ in $\C^2$. Then $F$
holomorphically extends to a neighborhood
of the whole curve $S$ in $\C^2$.

Proof.}
Without loss of generality, in a local
coordinate system $z=(z_1, z_2)=x+iy$,
$M$ and $S$ have the equations
$x_2=0$ and $z_2=0$ respectively,
and $p=(0,0)$.
Choose $r_1=r_2$ and $E_2$ so small that
$F$ is holomorphic in a neighborhood of
$I_1\times E_2$. Then by Lemma 5, $F$ extends
to a neighborhood  of $I_2\times E_1$ for
big $E_1$. The proposition is proved.
\medskip

{\bf Theorem 6.
\sl
Let $\L\subset \R^2$ be a real-analytic convex curve
with strictly positive curvature.
Let $f$ be a complex function in the exterior of $\L$.
Let $l_\l$ be the tangent line to $\L$ at $\l\in\L$.
Suppose for every $\l\in\L$, there is an entire function
$f_\l$ on $\C=\R^2$ such that $f|_{l_\l}=f_\l|_{l_\l}$.
Then $f$ extends as an entire function on $\C^2$.

Proof.}
We need to show that if $\L$ is real-analytic, then
the function $F$ in the proof of Theorem 1 is
continuous on $M$.
We first show that $F$ is real-analytic on
$\Pi_e\subset\Sigma$.
Following [\"O], we define
$$
\phi: Y=(\R\times\R)\setminus K\to\Pi_e,\quad
\phi(s_1,s_2)=h_{\l(s_1)}\cap h_{\l(s_2)},
$$
where $K$ is the set of all $(s_1,s_2)$ such that
the lines $h_{\l(s_1)}$ and $h_{\l(s_2)}$ are parallel.
Since $\L$ is real-analytic, then so is $\phi$.
Then $\phi$  extends to a holomorphic map $\Phi$ in
a neighborhood of $Y$ in $\C^2$.
Let $p\in\Pi_e$, $p=\phi(c_1, c_2)$.
We use $\Phi^{-1}$ to introduce local coordinates
$(s_1+it_1,s_2+it_2)$ in a neighborhood of $p$.
In these coordinates $M$ is represented as
the union of the hyperplanes $s_1=c_1$ and $s_2=c_2$,
on which $F$ is separately holomorphic.
By Lemma 4 the function $F$ holomorpically extends
to a neighborhood of $p$.
Finally, by Proposition 5, the analyticity of $F$
propagates along complex lines in $M$,
so $F$ is real-analytic, whence continuous on $M$.
The theorem is proved.

\beginsection References

\frenchspacing

\item{[AEK]}
V. Aguilar, L. Ehrenpreis, P. Kuchment,
Range condition for the exponential Radon transform,
{\sl J. d'Analyse Math. \vol 68}
(1996), 1--13.

\item{[K]}
H. Kneser,
Die Randwerte einer analytischen Funktion zweier Veränderlichen,
{\jour Monatsh. Math. Phys. \vol 43}
(1936), 364-–380.

\item{[L]}
H. Lewy,
On the local character of the solutions of an atypical
linear differential equation in three variables and a
related theorem for regular functions of two complex
variables,
{\jour Ann. Math. \vol 64}
(1956), 514--522.

\item{[\"{O}]}
O. \"{O}ktem,
Extension of separately analytic functions and applications
to range characterization of the exponential Radon transform,
{\sl Ann. Polon. Math. \vol 70}
(1998), 195--213.

\item{[S1]}
J. Siciak,
Analyticity and separate analyticity of functions defined
on lower dimensional subsets of $\C^n$,
{\sl Zeszyty Nauk. Uniw. Jagiello. Prace Mat. Zeszyt \vol 13}
(1969), 53--70.

\item{[S2]}
J. Siciak,
Separately analytic functions and envelopes of holomorphy
of some lower dimensional subsets of $C^n$,
{\sl Ann. Polon. Math. \vol 22}
(1969/1970), 145--171.

\item{[T]}
A. Tumanov,
Testing analyticity on circles,
{\sl Amer. J. Math. \vol 129}
(2007), 785--790.

\bigskip\bigskip\bigskip

\noindent
{\auth Alexander Tumanov}\br
Department of Mathematics, University of Illinois\br
1409 W. Green Street, Urbana, IL 61801, U. S. A.\br
{\sl E-mail: tumanov@uiuc.edu}

\bye